\documentclass[twocolumn]{autart}    
\usepackage{amssymb}
\usepackage[dvips]{epsfig}
\usepackage{amsmath,amssymb}
\usepackage{latexsym}
\usepackage{subfigure}
\usepackage{color}
\usepackage[center]{caption}
\usepackage{graphicx}
\usepackage{url}

\newcommand{\be}{\begin{equation}}
\newcommand{\eq}{\end{equation}}
\newcommand{\ba}{\begin{array}}
\newcommand{\ea}{\end{array}}
\newcommand{\bean}{\begin{eqnarray*}}
\newcommand{\eean}{\end{eqnarray*}}
\newcommand{\bea}{\begin{eqnarray}}
\newcommand{\eea}{\end{eqnarray}}

\newcommand{\R}{\rm {I\kern-2pt R}}

\newcommand{\beq}{\begin{equation}}
\newcommand{\eeq}{\end{equation}}

\newtheorem{theorem}{\bf Theorem}[section]
\newtheorem{remark}{\bf Remark}[section]
\newtheorem{lemma}{\bf Lemma}[section]
\newtheorem{example}{\bf Example}[section]

\newtheorem{corollary}{\bf Corollary}[section]

\newtheorem{definition}{\bf Definition}[section]

\begin{document}

\begin{frontmatter}

\title{Synchronization of High-Dimensional Linear Networks over Finite Fields} 
\thanks[footnoteinfo]{This work is supported in part by National Natural Science Foundation
(NNSF) of China under Grants 62373193 and 62103194, China Postdoctoral Science Foundation under Grant 2023M731692, and Postgraduate Research $\&$ Practice Innovation Program of Jiangsu Province of China.
}

\author[]{Siyu Zou}\ead{symathz@163.com},  
\author[]{Ting Li}\ead{dylliting@njnu.edu.cn},  %
\author[]{Jiandong Zhu*}\ead{zhujiandong@njnu.edu.cn}    %

\address{Key Laboratory of NSLSCS (Ministry of Education), School of Mathematical Sciences, Nanjing Normal University, Nanjing,
        210023, PRC}  
%
%
%
%
%

\begin{keyword}                           
Synchronization; Consensus; High-dimensional networks; Finite fields          
\end{keyword}                             
\begin{abstract}
	This paper investigates the synchronization problems for general high-dimensional linear networks over finite fields. By using the technique of linear transformations and invariant subspaces for linear spaces over finite fields, several necessary and sufficient conditions for the synchronization of high-dimensional linear networks over finite fields are proposed. This paper not only generalizes the existing results from 1-dimensional to high-dimensional linear networks but also adopts a new approach. Finally, some numerical examples are given to illustrate the effectiveness of our theoretical results.
\end{abstract}
\end{frontmatter}

\section{Introduction}
Linear dynamical systems over finite fields have been studied for a long time as mathematical models of logical circuits, genetic networks and so on  \cite{Elspas1959}, \cite{LAUBENBACHER2004523}, \cite{Toledo2004}.
If a dynamical system is composed of many agents interconnected by a network, it is also called a dynamical network, whose dynamical properties are usually determined by the dynamical equations of the agents and the network topology. The consensus problem for linear dynamical networks over finite fields is investigated in \cite{consensus_Pasqualetti}, which has important applications in distributed task privacy for aggregation \cite{Distributed-Connor} and
distributed computing over encrypted data
\cite{Distributed-Freris}. The so called consensus means that all the agent states enter the same fixed point, which is possibly dependent on the initial states.
In \cite{consensus_Pasqualetti}, a main contribution is that a necessary and sufficient condition described by the characteristic polynomial for consensus is proposed.
Synchronization, a generalization of consensus, for linear networks over finite fields is studied in \cite{synchronization4}. Different from consensus, the concept of synchronization only requires that the state errors among the agents enter the origin, which needs the terminal state to be a fixed point.
In \cite{synchronization4}, as a generalization of \cite{consensus_Pasqualetti}, a similar necessary and sufficient condition described by the characteristic polynomial for synchronization is derived. \cite{Wang-Synchronous}, an equivalent condition to synchronization is proposed based on the concept of inverse recursion subspace. Further generalized contributions on linear networks over finite fields have been achieved. For example, switching topologies or time delays are considered in \cite{consensus3}, \cite{STP1}, \cite{STP2}, \cite{synchronization6},
cluster synchronization is investigated in
 \cite{synchronization5}, and stochastic networks over finite fields are studied in \cite{optimal}, \cite{Lin-Synchronization},\cite{Lin-Quotients}.

The linear networks addressed in the above literature are composed of one-dimensional agents, and thus are usually called one-dimensional linear networks. If a network's every agent has an $m$-dimensional state, we call it an $m$-dimensional network.
In \cite{leaderconsensus}, \cite{synchronization4}, \cite{Lin-Leader-Follower} and   \cite{Synchronizability_Wang}, high-dimensional linear networks over finite fields have been studied. But the network topologies and the system model have special forms. To the best of our knowledge, the consensus and synchronization problems for the general form of high-dimensional linear networks over finite fields have not been addressed yet.

In this paper, we investigate the synchronization and consensus problems for general high-dimensional linear networks over finite fields. A necessary and sufficient condition for synchronization described by the characteristic polynomial is proposed. A simple linear transformation and an invariant subspace technique are introduced to deal with the high-dimensional synchronization problem. Compared with the existing literature on 1-dimensional linear networks, our method is simpler and straightforward. The constructed invariant subspace plays an important role in overcoming the essential difficulty for the synchronization of the general high-dimensional linear networks. In particular, the problem of consensus is accordingly solved as an application of the obtained synchronization result. Finally, some numerical examples are provided to illustrate the main contributions.

The rest of this paper is organized as follows; Section 2 contains some preliminaries and the problem statement. In Section 3, the main results are presented. Section 4 includes some numerical examples. Section 5 is a brief conclusion.

\section{Preliminaries and problem statement}
In this section, we give some necessary preliminaries and problem statement. Firstly, we present the following notations.

\begin{itemize}
	
	\item $\mathbb{N}$: the set of positive integers;
	
	\item $\mathbb{F}_p$: finite field with $p$ elements and $p$ is a prime number;
	
	\item $\mathbb{F}_p^{n}$: the set of $n$-dimensional column vectors over finite field $\mathbb{F}_p$;
	
	\item $\mathbb{F}_p^{n\times n}$: the set of $n \times n$ matrices over finite field $\mathbb{F}_p$;
	
	\item $\mathbf{1}_n$: $n$-dimensional column vector of all ones;
	
	\item $I_{n}$: $n$-dimensional identity matrix;
	
	\item $\dim (\mathbb{W})$: the dimension of linear space $\mathbb{W}$;
	
	\item $\deg(f)$: the degree of the polynomial $f$;
	
	\item $\otimes$: the Kronecker product.
\end{itemize}

Consider a network with $n$ nodes, whose states are defined on finite field $\mathbb{F}_p$ with $p$ being a prime number. Each node $i$ represents an agent:
\begin{equation}
	\label{xi}
	x_i(t+1)=A_{i1}x_1(t)+A_{i2}x_2(t)+ \cdots +A_{in}x_n(t),
\end{equation}
where $x_i(t)=[x_{i1}(t),\cdots,x_{im}(t)]^{\mathrm{T}} \in \mathbb{F}_p^m$ is the state of the agent $i$ at time $t$,  $A_{ij} \in \mathbb{F}_p^{m \times m}$ is the coupling matrix for each $i,j=1,2,\cdots n$.

Thus the dynamics of the $m$-dimensional linear network over finite field $\mathbb{F}_p$ is described by
\begin{equation}
	\label{mymodel}
	x(t+1)=Ax(t),
\end{equation}
where  $x(t)=[x_1^{\mathrm{T}}(t),x_2^{\mathrm{T}}(t),\cdots,x_n^{\mathrm{T}}(t)]^{\mathrm{T}} \in \mathbb{F}_p^{nm}$ and $A=(A_{ij}) \in \mathbb{F}_p^{nm \times nm}$.

\begin{remark}
In \cite{synchronization4} and \cite{synchronization6}, the considered $m$-dimensional linear network over finite field has a special form as $x(t+1)=(M \otimes I_m)x(t)$. For this particular case, the results on  $1$-dimensional linear network can be easily generalized. But for general case \eqref{mymodel}, the generalization is not trivial.
\end{remark}

Similar to the definitions in \cite{synchronization4} and \cite{consensus_Pasqualetti}, the definitions of synchronization and consensus of high-dimensional linear network \eqref{mymodel} over finite field $\mathbb{F}_p$ are stated as follows:

\begin{definition}
	\label{def1}
Linear network \eqref{mymodel} over $\mathbb{F}_p$ is said to achieve synchronization if for any initial state vector $x(0) \in \mathbb{F}_p^{nm}$, there exists a finite time $T \in \mathbb{N}$ such that $x_1(t)=x_2(t)= \cdots =x_n(t)$ for all $t \ge T$.
\end{definition}

\begin{definition}
	\label{def2}
 Linear network \eqref{mymodel} over $\mathbb{F}_p$ is said to achieve consensus if for any initial state vector $ x(0) \in \mathbb{F}_p^{nm}$, there exists a finite time $T \in \mathbb{N}$ such that $x(T)=x(K+\tau)=\mathbf{1}_n \otimes \alpha$ for all $\tau \in \mathbb{N} $ and for some $\alpha \in \mathbb{F}_p^m$.
\end{definition}

Obviously, a consensus network is a special synchronized network.
The state vector of a consensus network reaches a fixed point, but the state vector of a synchronized network probably enters a cycle.

The \emph{transition graph} associated with network \eqref{mymodel} over $\mathbb{F}_p$ is defined as $\mathcal{G}_A=(\mathcal{V}_A,\mathcal{E}_A)$ with the vertex set $\mathcal{V}_A=\left\{ v \mid v\in \mathbb{F}_p^{nm}\right\}$ and the edge set $\mathcal{E}_A=\left\{(u,v) \mid Av=u, \forall u,v \in \mathcal{V}_A \right\}$. Essentially, the transition graph $\mathcal{G}_A$ represents the trajectory of network $x(t+1)=Ax(t)$. Denote the set of vertices in all cycles in $\mathcal{G}_A$ as $\mathbb{C}_A$. Let $P_A(\lambda)$ be the characteristic polynomial of $A$ and it can be factorized as $P_A(\lambda)=\lambda^kf_A(\lambda)$ with $f_A(0) \neq 0$. It has been shown in \cite{Toledo2004} that $\mathbb{C}_A$ is a linear subspace of $\mathbb{F}_p^{nm}$, and the trajectory of network \eqref{mymodel} will enter $\mathbb{C}_A$ in $nm$ steps at most. From \cite{Toledo2004}, the following lemmas follows:

\begin{lemma}
	\label{lemma-dim}
	$\dim(\mathbb{C}_A)=\deg[f_A(\lambda)]$.
\end{lemma}

\begin{definition}
	A matrix $M \in \mathbb{F}_p^{r \times r}$ for any $r \in \mathbb{N}$ is called a nilpotent matrix if there exists a positive integer $k$ such that $M^k=0_{r \times r}$.
\end{definition}

\begin{lemma}
	\label{lem4}
Consider the dynamical equation
\begin{equation}
\label{eq-w}
 w(t+1)=Mw(t), \ w(0)=w_0
\end{equation}
over a finite field $\mathbb{F}_p$, where $w\in \mathbb{F}_p^r$ and $M\in \mathbb{F}_p^{r\times r}$. The following statements are equivalent:
\\
\mbox{\rm~~(i)} the state $w(t)$  terminates at $0_r$ for every initial condition $w(0)$;\\
\mbox{\rm~~(ii)} the characteristic polynomial $P_M(\lambda)$ is $\lambda^r$; \\
\mbox{\rm~~(iii)} the coefficient matrix $M$ is a nilpotent matrix.
\end{lemma}
%
%
%
%

\begin{lemma}
	\label{lem5}
Consider the linear system \eqref{eq-w} over a finite field $\mathbb{F}_p$. For each initial condition $w(0)=w_0\in \mathbb{F}_p^r$, there is a constant vector $\alpha(w_0)$ such that $w(t)$ terminates at $\alpha(w_0)$ if and only if the minimum polynomial of $M$ is in the form of $m_M(\lambda)=\lambda^{s}(\lambda-1)$ for some $0\leq s\leq r-1$.
\end{lemma}

\section{Main Results }
For studying synchronization of network \eqref{mymodel} over $\mathbb{F}_p$, we define a synchronization set as follows:
 \begin{equation}
\label{eq1}
\mathbb{S}:=\left\{\mathbf{1}_n \otimes \alpha \mid \alpha \in \mathbb{F}_p^m\right\}.
\end{equation}
By Definitions \ref{def1} and \ref{def2}, one can easily obtain the lemma as follows:
\begin{lemma}
	\label{lem2}
	For the network \eqref{mymodel} over the finite field $\mathbb{F}_p$, we have
\\ \mbox{\rm~~(i)} The network \eqref{mymodel} over $\mathbb{F}_p$ achieves synchronization if and only if $\mathbb{C}_A \subseteq \mathbb{S}$. \\ \mbox{\rm~~(ii)} The network \eqref{mymodel} over $\mathbb{F}_p$ achieves consensus if and only if $\mathbb{C}_A \subseteq \mathbb{S}$ and $\mathbb{C}_A=\left\{ \beta \in \mathbb{F}_p^{nm} | A \beta=\beta \right\}$.
\end{lemma}
\subsection{The case that the synchronization set is $A$-invariant}
Consider the network \eqref{mymodel} over $\mathbb{F}_p$ in the special case that the synchronization set $\mathbb{S}$ is $A$-invariant, i.e. $A\mathbb{S}\subseteq\mathbb{S}$. Before the main result of this subsection, we first give a lemma:
\begin{lemma}
\label{lem3.1}
Let $A_i=\sum_{j=1}^nA_{ij}$. Then the synchronization set $\mathbb{S}$ is $A$-invariant if and only if $A_1=A_2=\cdots=A_n$, i.e. $A(\mathbf{1}_n\otimes I_m)=\mathbf{1}_n\otimes A_1$.
\end{lemma}
\vspace{-4mm}
\begin{pf}
Since $\mathbb{S}$ has the form of (\ref{eq1}), we have that $A\mathbb{S}\subseteq \mathbb{S}$ if and only if, for any $\alpha\in \mathbb{F}_p^m$, there exists $\beta \in \mathbb{F}_p^m$ such that
$
A(\mathbf{1}_n\otimes \alpha)=\mathbf{1}_n\otimes \beta,
$
that is,
\begin{equation}
\label{eq}
\begin{bmatrix}
A_1 \alpha\\[-2pt]
 A_2 \alpha	\\[-2pt]
  \vdots
	\\[-2pt]
 A_n \alpha
		\end{bmatrix} =
\begin{bmatrix}
		\beta
			\\[-2pt] \beta
			\\[-2pt] \vdots
			\\[-2pt] \beta
		\end{bmatrix}.
\end{equation}
Due to the arbitrarity of $\alpha\in \mathbb{F}_p^m$, we conclude that (\ref{eq}) holds if and only if $A_1=A_2=\cdots=A_n$.
	\hfill $\blacksquare$
\end{pf}

\begin{theorem}
	\label{th1}
	Consider the network \eqref{mymodel} over the finite field $\mathbb{F}_p$. If the synchronization set $\mathbb{S}$ described by (\ref{eq1}) is $A$-invariant, then network \eqref{mymodel} achieves synchronization if and only if $P_A(\lambda)=\lambda^{nm-m}P_{A_1}(\lambda)$.
\end{theorem}
\vspace{-3mm}
\begin{pf}
For studying synchronization, we construct a nonsingular transformation for  \eqref{mymodel} as follows:
$$ z_1(t)=x_1(t),\ \ z_i(t)=x_i(t)-x_1(t), i=2,3,\dots,n,
$$
which can be expressed in a compact form $z(t)=Tx(t)$ with
$$
T\!\!=\!\!\begin{bmatrix}
	I_m &  0 \\
	\!-\mathbf{1}_{n-1} \!\!\otimes \!\!I_m & I_{(\!n-1\!)m}  \end{bmatrix},\ \
T^{-1}\!\!=\!\!\begin{bmatrix}
I_m &  0 \\
\!\mathbf{1}_{n-1} \!\!\otimes \!\!I_m & I_{(\!n-1\!)m}  \\
\end{bmatrix}.
$$
Considering the structure of the partitioned matrices, we write $A$ in the form as follows:
$$
A=\begin{bmatrix}
A_{11} &  \tilde A_{12} \\
\tilde A_{21} & \tilde A_{22} \\
\end{bmatrix}.
$$
Then network \eqref{mymodel} can be transformed into
\begin{equation}
	\label{eq5}
	z(t+1)=TAT^{-1}z(t).
\end{equation}
Since $\mathbb{S}$ is $A$-invariant, by Lemma
\ref{lem3.1}, we have $$A(\mathbf{1}_n\otimes I_m)=\mathbf{1}_n\otimes A_1.$$
Thus, a straightforward computation shows that
 \begin{eqnarray}
	\label{eq6}
T\!AT^{-\!1}&\!\!=&\!\!\begin{bmatrix}
	I_m &  0 \\
	\!-\!\mathbf{1}_{\!n-1} \!\!\otimes \!\!I_m & I_{(\!n-1\!)m}\!  \end{bmatrix}
\!\!\!\begin{bmatrix}
\!A_{11} &  \tilde A_{12}\! \\
\!\tilde A_{21} & \tilde A_{22}\!
\end{bmatrix}	
\!\!\!\begin{bmatrix}
	I_m &  0 \\
	\!\mathbf{1}_{\!n-1} \!\!\otimes \!\!I_m & I_{(\!n-1\!)m}\!  \end{bmatrix}\nonumber \\
&\!\!=&\!\!\begin{bmatrix}
	I_m &  0 \\
	\!-\!\mathbf{1}_{\!n-1} \!\!\otimes \!\!I_m & I_{(\!n-1\!)m}\!  \end{bmatrix}
\!\!\begin{bmatrix}
\!A_{1} &  \tilde A_{12}\! \\
\mathbf{1}_{\!n-1}\!\!\otimes \!\! A_1 & \tilde A_{22}\!
\end{bmatrix}	
\nonumber \\
&\!\!=&\!\!\begin{bmatrix}
	A_{1} & \tilde A_{12}  \\
	0 & \tilde A_{22}\!-\!\mathbf{1}_{\!n-1}\!\otimes \!\tilde A_{12}  \end{bmatrix}.
\end{eqnarray}
Therefore, by \eqref{eq6}, network \eqref{eq5} can be rewritten as
\begin{equation}
		\left\{\begin{array}{l}
z_1(t+1)=A_1z_1(t)+ \tilde A_{12} \tilde{z}(t) \\	\tilde{z}(t+1)=\bar A_{22} \tilde{z}(t),
		\end{array}\right.,	
	\end{equation}
where $\bar A_{22}=\tilde A_{22}-\mathbf{1}_{\!n-1}\!\otimes \!\tilde A_{12}$,  $\tilde{z}(t)=[z_2^{\mathrm{T}}(t),\dots,z_n^{\mathrm{T}}(t)]^{\mathrm{T}} \in \mathbb{F}_p^{(n-1)m}$.

Based on the above analysis, we conclude that network \eqref{mymodel} achieves synchronization if and only if every trajectory of the $\tilde{z}$-subnetwork terminates at $0_{(n-1)m}$, in other words, the matrix $\bar A_{22}$ is nilpotent over the finite field $\mathbb{F}_p$, i.e. $P_{\bar A_{22}}(\lambda)=\lambda^{nm-m}$, which is equivalent to
$$
P_A(\lambda)\!=\!P_{TAT^{-1}}(\lambda)\!=\!P_{\bar A_{22}}(\lambda)P_{A_1}(\lambda)\!=\!\lambda^{nm-m}P_{A_1}(\lambda).
$$
	\hfill $\blacksquare$
\end{pf}
\vspace{-4mm}
For the special case of $m=1$, by Theorem \ref{th1} and Lemma \ref{lem3.1}, we obtain the following Corollary \ref{cor1},
which is just Theorem 3 of \cite{synchronization4}:
\begin{cor}
	\label{cor1}
Consider the network \eqref{mymodel} over the finite field $\mathbb{F}_p$ with $m=1$. If there exists a scalar $\alpha$ such that $A\mathbf{1}_n=\alpha\mathbf{1}_n$, then network \eqref{mymodel} achieves synchronization if and only if $P_A(\lambda)=\lambda^{n-1}(\lambda-\alpha)$.
\end{cor}

\begin{cor}
	\label{cor2}
Assume that network \eqref{mymodel} satisfies the conditions of Theorem \ref{th1}. Then the final synchronization behavior is determined by
$x_1(t+1)=A_1x_1(t)$. In particular, if network \eqref{mymodel} has an $A$-invariant synchronization set, then it achieves consensus if and only if $P_A(\lambda)=\lambda^{nm-m}P_{A_1}(\lambda)$, and $A_1$ has a minimum polynomial $\lambda^{s}(\lambda-1)$
for some $0\leq s\leq m-1$.
\end{cor}

\begin{remark}
 Our method is completely different from and much simpler than that in \cite{consensus_Pasqualetti} and \cite{synchronization4}. In the existing literature, the cycle structure of linear systems over finite fields described by Theorem 5 of \cite{Toledo2004} plays an important role in the necessity proof. However, we only use a simple linear transformation instead of the complex cycle structure to solve the synchronization problem.
\end{remark}
This subsection is discussed under the condition that $\mathbb{S}$ is $A$-invariant.
Actually, for the network \eqref{mymodel} with $m=1$ and a non-nilpotent $A$, the invariance of $\mathbb{S}$ with respect to $A$ is also necessary for synchronization. However, for high-dimensional linear networks, it is not necessary to require that $\mathbb{S}$ is $A$-invariant (see Example \ref{exa2}), which is an essential difference between one-dimensional networks and high-dimensional networks. Therefore, in the next subsection, we focus on the general case without the invariance condition of the synchronization set.

\subsection{The general case}
To explore synchronization conditions for the general high-dimensional linear networks over  finite fields, we construct two linear subspaces for the linear network \eqref{mymodel} over the finite field $\mathbb{F}_p$ as follows:
\begin{eqnarray}
\label{eqq9}
	\ \ \ \mathbb{W}_1&:=&\{\alpha\in \mathbb{F}_p^m \mid\ (A_i-A_1)A_1^t\alpha=0_m, \ i=2,\dots,n, \nonumber \\
&&\hspace{12mm}\ \  \ \ t =0,1,\dots,m\!-\!1\}. \\
\label{eq9}
	 \mathbb{W}&:=&\{ \mathbf{1}_n \otimes \alpha \mid \alpha\in \mathbb{W}_1\},
\end{eqnarray}
where $A_i=A_{i1}+A_{i2}+\cdots+A_{in} \in \mathbb{F}_p^{m \times m}$ for all $i=1,2,\dots,n$.

\begin{lemma}
\label{lem3.2}
Linear subspace $\mathbb{W}_1$ can be rewritten as
\begin{eqnarray}
\label{eq9}
	\ \ \ \mathbb{W}_1&:=&\{\alpha \in \mathbb{F}_p^m\mid A_1^t\alpha=A_2^t\alpha=\cdots=A_n^t\alpha, \  \nonumber \\
&&\hspace{16mm} t =1,2,\dots,m\}.
\end{eqnarray}
\end{lemma}

\begin{lemma}
\label{lem3.3}
$\mathbb{W}_1$ is $A_1$-invariant, and	$\mathbb{W}$ is $A$-invariant, i.e. $A_1\mathbb{W}_1 \subseteq \mathbb{W}_1$ and
$A\mathbb{W} \subseteq \mathbb{W}$.
\end{lemma}
\vspace{-5mm}
\begin{pf}
By Cayley-Hamilton Theorem, we have that $A_1^m$ can be linearly represented by $I_m$, $A_1$, $A_1^2$, $\dots,$ $A_1^{m-1}$. Thus, it follows from \eqref{eqq9} that $\mathbb{W}_1$ is $A_1$-invariant.
For any $\mathbf{1}_n\otimes \alpha\in \mathbb{W}$, we have $ A_1\alpha=A_2\alpha=\cdots=A_n\alpha$.
Hence, a straightforward computation shows that
\begin{equation}
\label{eqq12}
A(\mathbf{1}_n\otimes \alpha)=
\left[\begin{matrix}
  A_1\alpha \\[-2pt]
  A_2\alpha \\[-2pt]
  \vdots\\[-2pt]
  A_n\alpha
\end{matrix}\right]=\mathbf{1}_n\otimes A_1\alpha.
\end{equation}
From $A_1\mathbb{W}_1 \subseteq \mathbb{W}_1$, it follows that $A_1\alpha\in \mathbb{W}_1$.
Therefore, by (\ref{eqq12}), we conclude that
$A(\mathbf{1}_n\otimes \alpha)\in \mathbb{W}$, which means that $\mathbb{W}$ is $A$-invariant.
	\hfill $\blacksquare$
\end{pf}

\begin{lemma}
 	\label{lem3.4}
The linear network \eqref{mymodel} achieves synchronization if and only if the trajectory of \eqref{mymodel} for any initial state enters $\mathbb{W}$.
\end{lemma}
\vspace{-3mm}
\begin{pf}
Since $\mathbb{W}\subseteq \mathbb{S}$, the sufficiency is obvious. Now, we prove the necessity. Since the network achieves synchronization, every trajectory of $x(t)$ enters a cycle $\mathbb{C}\subseteq \mathbb{S}$. Letting $l=|\mathbb{C}|$, for any $\mathbf{1}_n\otimes \alpha\in \mathbb{C}$, one can write $\mathbb{C}$ as
$$
\mathbb{C}=\{A^i(\mathbf{1}_n\otimes \alpha)\mid \ i=0,1,2,\dots,l-1\}.
$$
Since $A(\mathbf{1}_n\otimes \alpha)\in \mathbb{S}$, there exists $\beta$ such that
$A(\mathbf{1}_n\otimes \alpha)=\mathbf{1}_n\otimes \beta$, which implies that
$$\beta=A_1\alpha=A_2\alpha=\cdots A_n\alpha.$$ Similarly, by $A^2(\mathbf{1}_n\otimes \alpha)\in \mathbb{S}$, we have $A(\mathbf{1}_n\otimes \beta)\in\mathbb{S}$, which implies that
$$A_1^2\alpha=A_2^2\alpha=\cdots A_n^2\alpha.$$
With the same arguments, we conclude that
$$A_1^i\alpha=A_2^i\alpha=\cdots A_n^i\alpha$$
for each $i \in \mathbb{N}$. Hence any $\mathbf{1}_n\otimes \alpha\in \mathbb{C}$ belongs to $\mathbb{W}$ due to Lemma \ref{lem3.2}.
Therefore, the trajectory of $x(t)$ for any initial state enters $\mathbb{W}$.
\hfill $\blacksquare$
\end{pf}

Let  $d=\dim(\mathbb{W})=\dim(\mathbb{W}_1)$ and denote a basis of $\mathbb{W}_1$ by
$\alpha_1,\alpha_2,\dots,\alpha_d$.
Then the set of vectors
$$
\beta_1\!:=\!\mathbf{1}_n\!\otimes\!\alpha_1, \ \beta_2\!:=\!\mathbf{1}_n\!\otimes\!\alpha_2,\  \dots, \ \beta_d\!:=\!\mathbf{1}_n\!\otimes\!\alpha_d
$$
is a basis of $\mathbb{W}$.
By Lemma \ref{lem3.3}, there exists $Q\in \mathbb{F}_p^{d\times d}$ such that
\begin{equation}
\label{eqq13}
A_1[\alpha_1,\alpha_2,\dots,\alpha_d]=[\alpha_1,\alpha_2,\dots,\alpha_d]Q.
\end{equation}
From \eqref{eqq12} and \eqref{eqq13}, it follows that
\begin{equation}
\label{eqq14}
A[\beta_1,\beta_2,\dots,\beta_d]=[\beta_1,\beta_2,\dots,\beta_d]Q.
\end{equation}
\begin{theorem}
	\label{th2}
	Consider the network \eqref{mymodel} over the finite field $\mathbb{F}_p$. The network \eqref{mymodel} over $\mathbb{F}_p$ achieves synchronization if and only if $P_A(\lambda)=\lambda^{nm-d}P_{Q}(\lambda)$, where $d=\dim(\mathbb{W})$ and $Q \in \mathbb{F}_p^{d \times d}$ determined by \eqref{eqq13} or \eqref{eqq14}.
\end{theorem}
\vspace{-3mm}
\begin{pf}
Extend the basis $\beta_1$, $\beta_2$, $\dots$, $\beta_d$ of the subspace $\mathbb{W}$ to a basis of $\mathbb{F}_p^{mn}$, and arrange these basis vectors into a matrix as follows:
\vspace{-2mm}
\begin{equation}
\label{eqq15}
T=[\beta_1,\ \beta_2, \dots, \beta_d,\ \beta_{d+1},\cdots,\ \beta_{mn}].
\end{equation}
Then, it follows from \eqref{eqq14} that
\begin{equation}
\label{eqq16}
AT=T\left[\begin{matrix}
      Q & \hat A_{12} \\[-2pt]
      0 & \hat A_{22}
    \end{matrix}\right],
\end{equation}
where $\hat A_{12}\in \mathbb{F}_p^{d\times (nm-d)}$ and $\hat A_{22}\in \mathbb{F}_p^{(nm-d)\times (nm-d)}$.
Take the nonsingular transformation $x=Tz$, which implies that $z$ is the coordinate vector of $x$ with respect to the basis shown in \eqref{eqq15}. Then, by \eqref{eqq16}, we get the dynamical equation of $z$ as follows:
\vspace{-2mm}
\begin{equation}
\label{eqq17}
\left[\begin{matrix}
     z_1(t+1)\\[-2pt]
    z_2(t+1)
    \end{matrix}\right]=	\left[\begin{matrix}
      Q & \hat A_{12} \\[-2pt]
      0 & \hat A_{22}
    \end{matrix}\right]	
    \left[\begin{matrix}
         z_1(t)\\[-2pt]
         z_2(t)
       \end{matrix}\right],	
	\end{equation}
where $z_1\in \mathbb{F}_p^d$ is composed of the coordinate components of $x$ with respect to the basis of $\mathbb{W}$, and $z_2\in \mathbb{F}_p^{nm-d}$ is composed of the other coordinate components. By Lemma \ref{lem3.4}, network \eqref{mymodel} achieves synchronization if and only if $z_2(t)$ terminates at $0_{mn-d}$, which is equivalent to that $\hat A_{22}$ is a nilpotent matrix. Further considering \eqref{eqq16}, we have that $\hat A_{22}$ is nilpotent if and only if $P_A(\lambda)=\lambda^{nm-d}P_{Q}(\lambda)$. Therefore, the proof is complete. \hfill $\blacksquare$
\end{pf}
\begin{corollary}
	\label{coro3.1}
	Consider the network \eqref{mymodel} over the finite field $\mathbb{F}_p$. The network \eqref{mymodel} over $\mathbb{F}_p$ achieves consensus if and only if $P_A(\lambda)=\lambda^{nm-d}P_Q(\lambda)$ and $Q$ has a minimum polynomial $\lambda^{s}(\lambda-1)$ for some $0\leq s\leq d-1$, where $d=\dim(\mathbb{W})$.
\end{corollary}
\vspace{-2mm}
\section{Numerical Examples}
This section gives three numerical examples to show the synchronization and consensus phenomena of high-dimensional networks over finite fields, respectively.
\begin{example}
	\label{exa1}
	\begin{figure}[t]
		\label{fig1}
		\centering
\includegraphics[height=6cm,width=6cm]{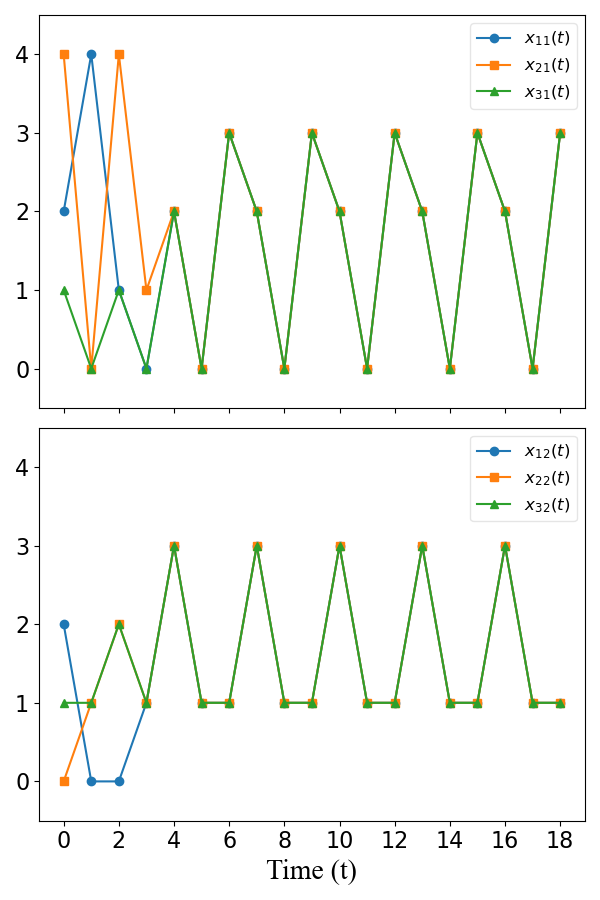}
		\caption{{The time response of the network  of Example \ref{exa1} when $x(0)=(2,2,4,0,1,1)^{\mathrm{T}}$. }}
		\label{fig11}
	\end{figure}
	Consider the 2-dimensional network  \eqref{mymodel} over $\mathbb{F}_5$ with 3 agents, where the coefficient matrix is
	\begin{equation}
		\label{eg-1}
		\begin{aligned}
			A=\begin{bmatrix}
				4 & 3 & 4 & 1 & 0 & 4 \\[-3pt]
				3 & 3 & 2 & 2 & 4 & 1 \\[-3pt]
				4 & 0 & 4 & 2 & 0 & 1 \\[-3pt]
				1 & 3 & 2 & 4 & 1 & 4 \\[-3pt]
				4 & 3 & 4 & 0 & 0 & 0 \\[-3pt]
				1 & 3 & 2 & 4 & 1 & 4
			\end{bmatrix},
		\end{aligned}
	\end{equation}
and the basic parameters are $m=2$, $n=3$ and  $p=5$. It is easily to obtain that
	\begin{equation}
		\notag
			A_1=A_2=A_3=\begin{bmatrix}
				3 & 3\\[-2pt]
				4 & 1
			\end{bmatrix}.
	\end{equation}
By Lemma \ref{lem3.1}, the synchronization set $\mathbb{S}$ of network \eqref{eg-1} is $A$-invariant. Since the characteristic polynomial of $A$ is $P_A(\lambda)=\lambda^4(\lambda^2+\lambda+1)$ and the characteristic polynomial of $A_1$ is $P_{A_1}(\lambda)=\lambda^2+\lambda+1$, the 2-dimensional network \eqref{eg-1} over $\mathbb{F}_5$ achieves synchronization by Theorem \ref{th1}. The time response sequences in Fig. \ref{fig11} show the synchronization.
\end{example}

\begin{example}
	\label{exa2}
	\begin{figure}[t]
		\label{fig2}
		\centering
		\includegraphics[height=6cm,width=6cm]{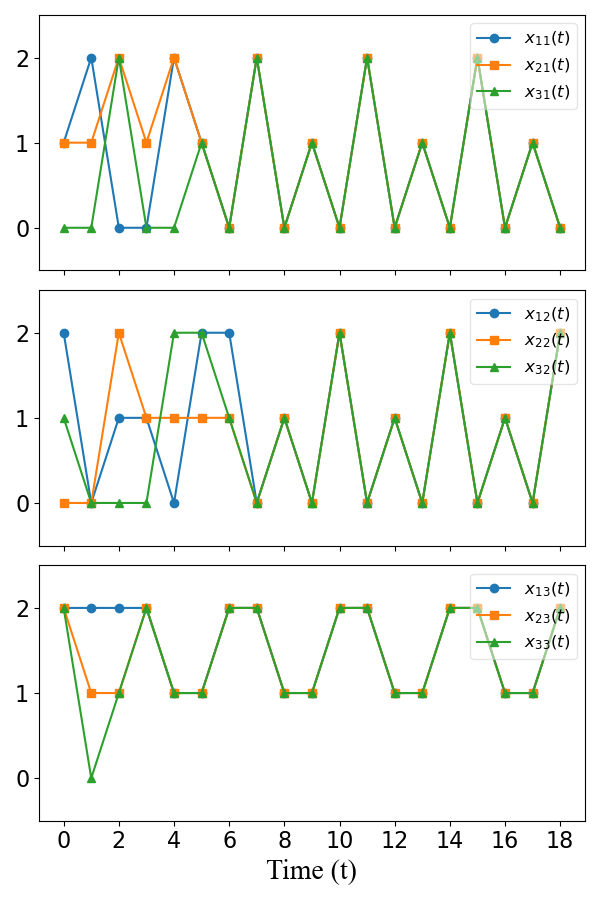}
		\caption{{The time response of the network in Example \ref{exa2} when $x(0)=(1,2,2,1,0,2,0,1,2)^{\mathrm{T}}$. }}
		\label{fig22}
	\end{figure}
Consider the 3-dimensional network  \eqref{mymodel} over $\mathbb{F}_3$ with $3$ agents and $m=3,n=3,p=3$, where the coefficient matrix is
	\begin{equation}
		\label{eg-2}
		\begin{aligned}
			A=\begin{bmatrix}
     0&0& 1& 0& 0& 1& 2& 0& 2\\[-2pt]
     0& 2& 0& 1& 1& 0& 0& 2& 1\\[-2pt]
     1& 0& 1& 2& 0& 2& 1& 0& 1\\[-2pt]
     2& 2& 0& 1& 1& 0& 1& 2& 2\\[-2pt]
     1& 0& 0& 0& 1& 0& 1& 2& 0\\[-2pt]
     0& 2& 2& 0& 1& 0& 0& 2& 0\\[-2pt]
     0& 1& 1& 0& 2& 0& 0& 1& 2\\[-2pt]
     0& 1& 0& 0& 0& 0& 0& 0& 2\\[-2pt]
     2& 1& 0& 1& 2& 2& 2& 1& 1			
\end{bmatrix}.
		\end{aligned}
	\end{equation}
With a straightforward calculation, we have
	\begin{equation}
		\notag
		\begin{aligned}
			A_1=\begin{bmatrix}
     2 & 0& 1\\[-2pt]
     1 &2 & 1\\[-2pt]
     1 & 0& 1
			\end{bmatrix},
			A_2=\begin{bmatrix}
     1  &   2 &    2\\[-2pt]
     2  &   0 &    0\\[-2pt]
     0  &   2 &    2
			\end{bmatrix},
			A_3=\begin{bmatrix}
	 0  &   1 &    0\\[-2pt]
     0  &   1 &    2\\[-2pt]
     2  &   1 &    0
    	\end{bmatrix}.
		\end{aligned}
	\end{equation}
By Lemma \ref{lem3.1}, the synchronization set $\mathbb{S}$ of network \eqref{eg-2} is not $A$-invariant. In this case, we need to calculate the
$A_1$-invariant set
$\mathbb{W}_1$ as follows:
	\begin{equation}
		\mathbb{W}_1=\ker\begin{bmatrix}
			A_2-A_1 \\[-2pt]
			A_3-A_1 \\[-2pt]
			(A_2-A_1)A_1 \\[-2pt]
			(A_3-A_1)A_1 \\[-2pt]
			(A_2-A_1)A_1^2 \\[-2pt]
			(A_3-A_1)A_1^2
		\end{bmatrix}
=\mathrm{span}
\begin{bmatrix}
	1 &1 \\
	2 &0\\
	0& 1
\end{bmatrix}.
	\end{equation}
Solving
    \begin{equation}
    	\notag
    	A_1\begin{bmatrix}
    	1 & 1 \\[-2pt]
    	2 &	0 \\[-2pt]
    	0 &	1
    	\end{bmatrix}=\begin{bmatrix}
    	1 & 1 \\[-2pt]
    	2 &	0 \\[-2pt]
    	0 &	1
    \end{bmatrix}Q,
    \end{equation}
we get $Q=\begin{bmatrix}
    	1 & 1\\[-2pt]
    	1 & 2
    \end{bmatrix} \in \mathbb{F} _3^{2 \times 2}.$
Since $P_A(\lambda)=\lambda^7(\lambda^2+1)$ and  $P_Q(\lambda)=\lambda^2+1$, by Theorem \ref{th2}, we conclude that the 3-dimensional network \eqref{eg-2} over $\mathbb{F}_3$ achieves synchronization. In Fig. \ref{fig22}, it shows that, from the initial state  $x(0)\!=\!(1,2,2,1,0,2,0,1,2)^{\!\mathrm{T}}$, state $x(t)$ enters the cycle:\\ $(2,0,2,2,0,2,2,0,2)^{\mathrm{T}} \rightarrow (0,1,1,0,1,1,0,1,1)^{\mathrm{T}}\rightarrow (1,0,1,1,0,1,1,0,1)^{\mathrm{T}} \rightarrow (0,2,2,0,2,2,0,2,2)^{\mathrm{T}} \rightarrow (2,0,2,2,0,2,2,0,2)^{\mathrm{T}}$.
\end{example}
\begin{example}
	\label{exa3}
	\begin{figure}[t]
		\label{fig3}
		\centering
		\includegraphics[height=6cm,width=6cm]{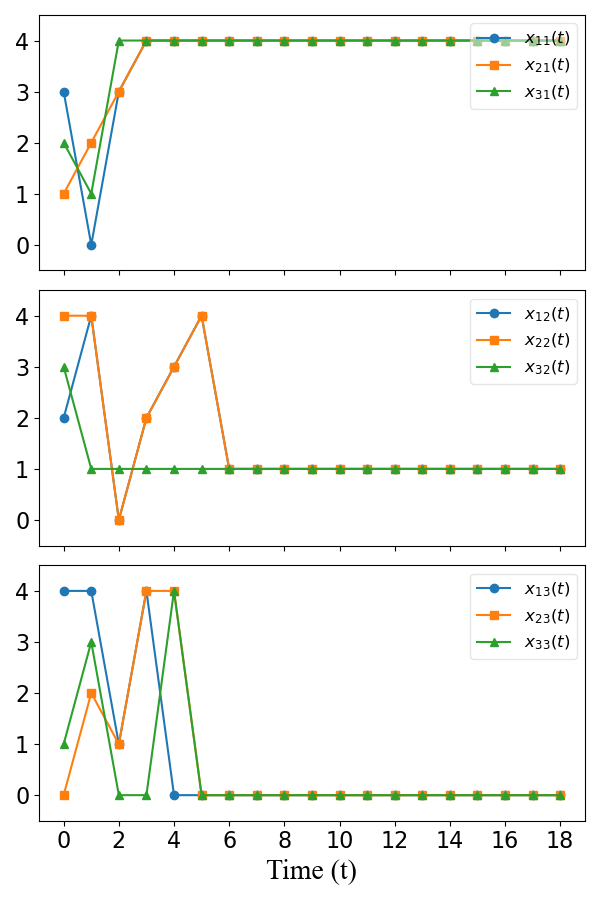}
		\caption{{The time response of the network of Example \ref{exa3} when $x(0)=(3,2, 4,1,4,0,2,3,1)^{\mathrm{T}}$. }}
		\label{fig33}
	\end{figure}
Consider the 3-dimensional network \eqref{mymodel} over $\mathbb{F}_5$ with 3 agents and $m=3,n=3,p=5$, where the coefficient matrix is
	\begin{equation}
		\label{eg-3}
		\begin{aligned}
			A=\begin{bmatrix}
 1 & 2& 0& 4& 3& 0& 3& 2& 0\\[-2pt]
 3& 1& 3& 4& 4& 0& 3& 1& 2\\[-2pt]
 0& 0& 1& 1& 0& 4& 4& 0& 1\\[-2pt]
 3& 3& 0& 2& 2& 0& 0& 4& 0\\[-2pt]
 2& 0& 3& 0& 0& 0& 2& 0& 2\\[-2pt]
 4& 4& 1& 3& 1& 0& 2& 4& 0\\[-2pt]
 2& 0& 0& 3& 0& 0& 2& 1& 0\\[-2pt]
 4& 0& 0& 1& 0& 0& 4& 0& 0\\[-2pt]
 2& 2& 1& 4& 3& 0& 1& 2& 0
			\end{bmatrix}.
		\end{aligned}
	\end{equation}
It is easily obtained that
		\begin{equation}
		\notag
		\begin{aligned}
			A_1=\begin{bmatrix}
				3	&2	&0 \\[-2pt]
				0	&1	&0 \\[-2pt]
				0	&0	&1
			\end{bmatrix},
			A_2=\begin{bmatrix}
				0	&4	&0 \\[-2pt]
				4	&0	&0 \\[-2pt]
				4	&4	&1
			\end{bmatrix},
			A_3=\begin{bmatrix}
				2	&1	&0 \\[-2pt]
				4	&0	&0 \\[-2pt]
				2	&2	&1
			\end{bmatrix}.
		\end{aligned}
	\end{equation}
With straightforward calculation, we get that
	\begin{equation}
	\mathbb{W}_1=\mathrm{span}
\begin{bmatrix}
	4 & 0\\[-2pt]
	1 & 0\\[-2pt]
	0 & 1
\end{bmatrix},\
Q=\begin{bmatrix}
    	1 & 0\\[-2pt]
    	0 & 1
    \end{bmatrix}.
    \end{equation}
Since the characteristic polynomial of $A$ is $P_A(\lambda)=\lambda^7(\lambda-1)^2=\lambda^7P_Q(\lambda)$ and the minimum polynomial of $Q$ is $\lambda-1$, the 3-dimensional network \eqref{eg-3} over $\mathbb{F}_5$ achieves consensus due to Corollary \ref{coro3.1}.
Fig. \ref{fig33} shows that, from the initial state $x(0)=(3, 2, 4, 1, 4, 0, 2, 3, 1)^{\mathrm{T}}$, state $x(t)$ terminates at $(4 ,1, 0, 4 ,1, 0, 4 ,1, 0)^{\mathrm{T}}$.
\end{example}

\section{Conclusion}
In this paper, we have investigated the synchronization problem for high-dimensional networks over finite fields. Necessary and sufficient conditions have been derived for synchronization and consensus, respectively. The obtained results have generalized the existing synchronization result in \cite{synchronization4} and a consensus result in \cite{consensus_Pasqualetti} for one-dimensional linear networks over finite fields. Our future research interests lie in exploring how to establish topology conditions for coupling networks, as well as designing synchronization protocols for individual systems.



\bibliographystyle{plain}        
\bibliography{zsy_refs}   

\begin{thebibliography}{10}

\bibitem{Elspas1959}
B.~Elspas.
\newblock The theory of autonomous linear sequential networks.
\newblock {\em IRE Transactions on Circuit Theory}, 6(1):45--60, 1959.

\bibitem{Distributed-Freris}
N.~M. Freris and P.~Patrinos.
\newblock Distributed computing over encrypted data.
\newblock In {\em 2016 54th Annual Allerton Conference on Communication,
  Control, and Computing (Allerton)}, pages 1116--1122, 2016.

\bibitem{LAUBENBACHER2004523}
R.~Laubenbacher and B.~Stigler.
\newblock A computational algebra approach to the reverse engineering of gene
  regulatory networks.
\newblock {\em Journal of Theoretical Biology}, 229(4):523--537, 2004.

\bibitem{consensus3}
X.~Li, M.~Z.~Q. Chen, H.~Su, and C.~Li.
\newblock Consensus networks with switching topology and time-delays over
  finite fields.
\newblock {\em Automatica}, 68:39--43, 2016.

\bibitem{STP1}
Y.~Li, H.~Li, and X.~Ding.
\newblock Set stability of switched delayed logical networks with application
  to finite-field consensus.
\newblock {\em Automatica}, 113:108768, 2020.

\bibitem{STP2}
Y.~Li, H.~Li, X.~Ding, and G.~Zhao.
\newblock Leader-follower consensus of multiagent systems with time delays over
  finite fields.
\newblock {\em IEEE Transactions on cybernetics}, 49(8):3203--3208, 2018.

\bibitem{optimal}
Y.~Li, H.~Li, and G.~Zhao.
\newblock Optimal state estimation for finite-field networks with stochastic
  disturbances.
\newblock {\em Neurocomputing}, 414:238--244, 2020.

\bibitem{Lin-Leader-Follower}
L.~Lin, J.~Cao, J.~Lam, S.~Zhu, S.~Azuma, and L.~Rutkowski.
\newblock Leader-follower consensus over finite fields.
\newblock {\em IEEE Transactions on Automatic Control}, page
  doi:10.1109/TAC.2024.3354195, 2024.

\bibitem{synchronization5}
L.~Lin, J.~Cao, X.~Liu, G.~Lu, and M.~Abdel-Aty.
\newblock Cluster synchronization of finite-field networks.
\newblock {\em IEEE Transactions on Cybernetics}, 2023.

\bibitem{Lin-Synchronization}
L.~Lin, J.~Cao, S.~Zhu, and P.~Shi.
\newblock Synchronization analysis for stochastic networks through finite
  fields.
\newblock {\em IEEE Transactions on Automatic Control}, 67(2):1016--1022, 2022.

\bibitem{Lin-Quotients}
L.~Lin, Z.~Jiang, H.~Lin, E.~C.~H. Ngai, and J.~Lam.
\newblock On quotients of stochastic networks over finite fields.
\newblock {\em IEEE Transactions on Control of Network Systems}, page
  doi:10.1109/TCNS.2023.3314583, 2024.

\bibitem{synchronization4}
M.~Meng, X.~Li, and G.~Xiao.
\newblock Synchronization of networks over finite fields.
\newblock {\em Automatica}, 115:108877, 2020.

\bibitem{Distributed-Connor}
M.~O'Connor and W.~B. Kleijn.
\newblock Distributed task privacy for aggregation using linear codes.
\newblock {\em IEEE Transactions on Signal and Information Processing over
  Networks}, 7:626--635, 2021.

\bibitem{consensus_Pasqualetti}
F.~Pasqualetti, D.~Borra, and F.~Bullo.
\newblock Consensus networks over finite fields.
\newblock {\em Automatica}, 50(2):349--358, 2014.

\bibitem{Toledo2004}
R.~A.~H. Toledo.
\newblock Linear finite dynamical systems.
\newblock {\em Communications in Algebra}, 33(9):2977--2989, 2005.

\bibitem{Synchronizability_Wang}
J.~Wang, J.~Feng, and Y.~Yu.
\newblock Synchronizability and protocol design of discrete-time second-order
  networks over finite fields.
\newblock {\em IEEE Transactions on Automation Science and Engineering}, page
  doi:10.1109/TASE.2024.3367886, 2024.

\bibitem{Wang-Synchronous}
J.~Wang, J.~Feng, Y.~Yu, and H.-L. Huang.
\newblock Synchronous networks over finite fields.
\newblock {\em IEEE Transactions on Automatic Control}, 68(11):6907--6912,
  2023.

\bibitem{leaderconsensus}
X.~Xu and Y.~Hong.
\newblock Leader-following consensus of multi-agent systems over finite fields.
\newblock In {\em 53rd IEEE Conference on Decision and Control}, pages
  2999--3004. IEEE, 2014.

\bibitem{synchronization6}
J.~Zhang, J.~Lu, M.~Xing, and J.~Liang.
\newblock Synchronization of finite field networks with switching multiple
  communication channels.
\newblock {\em IEEE Transactions on Network Science and Engineering},
  8(3):2160--2169, 2021.

\end{thebibliography}

\end{document}